\documentclass{IEEEtran}
\usepackage[utf8]{inputenc}

\usepackage{graphicx}
\usepackage{amsmath}
\usepackage{amssymb}
\usepackage{amsfonts}
\usepackage{mathrsfs}
\usepackage{dsfont}
\makeatletter
\def\amsbb{\use@mathgroup \M@U \symAMSb}
\makeatother

\usepackage[usenames,dvipsnames]{xcolor}
\usepackage{epstopdf}

\allowdisplaybreaks

\newtheorem{theorem}{Theorem}
\newtheorem{remark}{Remark}
\newtheorem{definition}{Definition}
\newtheorem{lemma}{Lemma}
\newtheorem{proposition}{Proposition}

\newtheorem{assumption}{Assumption}
\newtheorem{design condition}{Design condition}
\usepackage{float,graphicx,subfig,tikz}
\usetikzlibrary{shapes,arrows,scopes,plotmarks}

\newcommand{\vect}[1]{\boldsymbol{#1}} 
\usepackage{enumerate}

\usepackage{tabularx}
\DeclareMathOperator{\Ima}{Im}

\usepackage{balance}

\renewcommand{\baselinestretch}{1.08}

\begin{document}

\title{A distributed scheme for secondary frequency control with stability guarantees and optimal power allocation}

\author{Andreas Kasis\thanks{Andreas Kasis and Ioannis Lestas are with the Department of Engineering, University of Cambridge, Trumpington Street, Cambridge, CB2 1PZ, United Kingdom; e-mails: ak647@cam.ac.uk, icl20@cam.ac.uk},
Nima Monshizadeh\thanks{Nima Monshizadeh is with the Engineering and Technology Institute, University of Groningen, Nijenborgh 4, 9747AG, Groningen, The Netherlands. email: n.monshizadeh@rug.nl},
and Ioannis Lestas
\thanks{A preliminary version of this work has appeared in \cite{kasis2018ECC}.
 This manuscript includes the analytic proofs of the main results, additional  results,
  as well as further discussion, examples and simulations that demonstrate the significance of the proposed analysis. \newline
This work was supported by ERC starting grant 679774.}
}

\maketitle

\begin{abstract}
We consider the problem of distributed secondary frequency regulation in power networks such that stability and an optimal power allocation are guaranteed. 
This is a problem that has been widely studied in the literature, where two main control schemes have been proposed, usually referred to as {'Primal-Dual'} and
'distributed averaging proportional-integral (DAPI)' respectively.
 However, each has its limitations, with the former
{{incorporating additional} information flow requirements which may limit {its} applicability,}
  and {with} the existing literature on the latter {{relying} on static models for} generation and demand, {which is restrictive}.
We propose a novel control scheme that {aims to} overcome these issues by making use of generation measurements in {the control policy.} In particular, our controller relies on practical measurements and
allows distributed stability and optimality guarantees to be deduced for a broad range of
linear generation dynamics, {that can be of} higher order.
We show how the controller parameters can be selected in a computationally efficient way by solving appropriate linear matrix inequalities (LMIs).
 Furthermore, we demonstrate how the proposed analysis applies to various examples of turbine governor dynamics by using realistic numerical data.
 The practicality of our analysis is demonstrated with numerical simulations on the Northeast Power Coordinating Council (NPCC) 140-bus system that verify that our proposed controller achieves convergence to the nominal frequency, an economically optimal power allocation, and improved performance compared to existing schemes used in the literature.
\end{abstract}


\section{Introduction}\label{sec:Introduction}

\textbf{Motivation:}
{Current environmental concerns are {drawing} increasing attention on renewable sources of generation,}
with their penetration in power {networks expected to} grow over the next years \cite{lund2006large, ipakchi2009grid}.
The above will dramatically increase the number of active elements in the power network, making  its electromechanical behaviour less predictable and traditionally implemented, centralized control approaches   expensive and inefficient.
This highlights the importance of investigating distributed control schemes that will guarantee power network stability when such devices are included. {These concerns  have motivated recent studies} on distributed schemes with applications on both primary \cite{kasis2016primary, molina2011decentralized, kasis2019primary} and secondary frequency regulation \cite{trip2016internal, mallada2017optimal, kasis2019secondary}.

The introduction of highly distributed schemes for frequency {regulation} raises  an issue of economic optimality in the power allocation. Attempts to resolve this issue in {the literature resulted in} devising appropriately constructed optimization problems that ensured economic optimality and {designing the system equilibria in order} to be solutions to these problems. {It is evident in {the} literature that a synchronizing variable is useful for optimality to be achieved.
 While  frequency is used as the synchronizing signal in primary control studies \cite{kasis2016primary, zhaotopcu, devane2016primary}, some other signal, resulting from a suitably designed controller, {has been employed} for secondary frequency control (e.g. \cite{trip2016internal, kasis2017stability, dorfler}).}
However, when distributed optimal secondary frequency regulation is desired, the interaction of the physical system with the imposed communication scheme {may compromise the stability of the closed-loop system}. To cope with this, existing studies had to rely on restrictive assumptions, requiring generation to instantly follow a reference signal, or {measurement requirements that are hard to obtain,} such as continuous knowledge of demand,
  limiting the implementability of the proposed schemes.
Therefore, an open problem, which this study aims to address, is to obtain distributed stability and economic optimality guarantees for secondary frequency regulation, applicable to general network topologies,  without relying on restrictive measurement requirements.

\textbf{Literature survey:}
There are many recent studies associated with  {stability and optimality in distributed} secondary frequency control.
A common approach is to  involve control schemes with dynamics that follow from a primal/dual algorithm associated with some optimal power allocation optimization problem \cite{mallada2017optimal, kasis2017stability, low2014distributed, zhao2015distributed, li2016connecting, chen2018distributed, stegink2016unifying, wang2017distributed, zhao2018distributed, wang2018distributed}.
This approach allows to take into account economic considerations along {with} the objectives of secondary frequency control.
 {Furthermore, as has been demonstrated in \cite{kasis2017stability},} it allows for stability {and optimality} guarantees when high order and {nonlinear} generation dynamics {and convex cost functions} are considered.
  However, such schemes require knowledge of demand {in real time}, which can in some cases limit their practicality.
Attempts to adapt the Primal-Dual control scheme to avoid demand measurements
rely on additional information requirements \cite{li2016connecting, chen2018distributed, stegink2016unifying, wang2017distributed}, such as frequency derivative, system's inertia, damping coefficient and power transfers which may also lead to implementation challenges.

An alternative approach for optimal distributed secondary frequency regulation involves the use of distributed averaging proportional integral (DAPI) controllers \cite{trip2016internal, dorfler, dorfler2014plug, simpson2013synchronization, andreasson2014distributed}.
DAPI controllers are simple to implement,
requiring only knowledge of local frequency and exchange a synchronization signal without requiring any generation or load measurements.
On the other hand, existing results in {the literature incorporating DAPI controllers} do not accommodate
high order generation dynamics and restrict the {stability and optimality} analysis to static generation {and quadratic cost functions}.

{An alternative approach has been followed in \cite{trip2017distributed}, which proposes a controller that allows for stability {and optimality} guarantees when first and second order generation dynamics {and quadratic cost functions} are considered. However, the proposed scheme requires measurements of generators internal states which may be difficult to obtain. In addition, the imposed conditions in \cite{trip2017distributed} suggest conservative gains which may limit the performance and applicability of the results.
The authors in \cite{alghamdi2018conditions} consider  the effects of time delays in power networks with  second order turbine governor dynamics and consensus based distributed secondary frequency
control schemes inspired by \cite{kasis2018ECC, trip2017distributed}.}
{For a thorough survey of distributed approaches for stability and optimality in power systems, see \cite{molzahn2017survey}, \cite{dorfler2019distributed}.}

\textbf{Main contributions:}
In this paper, we propose a {{distributed control scheme for optimal secondary frequency regulation}, that will be referred to as} distributed averaging dynamic output control (DADOC). {A distinctive feature of this scheme is that it
{allows} for stability and optimality guarantees when high order generation dynamics are considered, {without} imposing {restrictive information flow requirements}.
{Our proposed scheme allows secondary frequency control to be performed in a
plug and play
 fashion and is applicable to general network topologies.}}

{DADOC schemes have} the advantage over DAPI schemes that they allow {the inclusion of} {higher order} generation dynamics, by imposing only an additional condition for knowledge of generation output.
{At the same time, they {impose less restrictive} measurement requirements,
 {compared to} Primal-Dual schemes that also provide stability guarantees when high order generation dynamics are considered.}
{Hence, DADOC controllers share advantages of both schemes,
allowing the inclusion of highly relevant generation dynamics with {easily obtainable measurement} requirements (generation and frequency).}

{Compared to the scheme proposed in \cite{trip2017distributed}, DADOC schemes allow for stability guarantees when {generation dynamics of arbitrary (finite) order are considered, while the work in \cite{trip2017distributed} is restricted to first and second order dynamics. Furthermore,
 second order generation dynamics
{with less restrictive stability conditions}
are included in the analysis
  and no internal state measurement requirements are imposed.}}

{Our analysis provides conditions for the design of the controller gains such that stability and optimality are guaranteed.}
An important feature {of} the proposed conditions, is that those can be  verified in a computationally efficient way by means of a linear matrix inequality (LMI). Several examples of relevant generation dynamics, {with realistic numerical data,} are provided to demonstrate the relevance of our {contribution.}

{Our analytic results are {accompanied by} numerical simulations on the NPCC 140-bus system
that demonstrate convergence to {the} nominal frequency and an economically optimal power allocation
at the presence of high order turbine governor dynamics.
 Furthermore, it is numerically demonstrated that DADOC schemes offer  improved performance compared to DAPI schemes, being able to provide significantly faster response to power disruptions.}

\textbf{Paper structure:} The  rest of the paper is structured as follows:
Section \ref{sec: Preliminaries} contains some basic preliminaries and in Section \ref{Problem_Formulation} we present the power network model and generation dynamics.
In Section \ref{sec:Power_Command_Dynamics} we present our proposed control scheme and conditions that allow an optimal power allocation at equilibrium.
 Section \ref{sec:Convergence} contains the main stability result of this paper, {which is associated with convergence to an optimal power allocation.}
In Section \ref{sec:Discussion} we clarify the importance of our proposed scheme compared to existing schemes in the literature and demonstrate its relevance with various applications on realistic generation models.
 Our results are validated {with numerical simulations} in Section \ref{Simulation_NPCC}.
 Finally, conclusions are drawn in Section \ref{Conclusions}.

\section{Preliminaries}\label{sec: Preliminaries}

 The set of $n$-dimensional vectors with real entries is denoted by $\mathbb{R}^n$.
The first derivative of a function $f(q)$, $f:\mathbb{R}\rightarrow \mathbb{R}$ is denoted by $f'(q) = \tfrac{d}{dq} f(q)$ {and its inverse by $f^{-1}(.)$.
A function $f:\mathbb{R}^n \rightarrow \mathbb{R}$ is said to be positive definite if $f(0) = 0$ and $f(x) > 0$ for every $x\neq0$.
 We write $\vect{0}_n$ and $\vect{1}_n$ to denote $n \times 1$ vectors with all elements equal to $0$ and $1$ respectively.
For a discrete set $S$, the term $|S|$ denotes its cardinality.
A matrix $A$ is said to be Hurwitz if all its eigenvalues lie on the open left half plane \cite{zhou1998essentials}. {A matrix $A$ is said to be positive definite (semi-definite) when $x^TAx > 0$ for all $x \neq 0$ (when $x^TAx \geq 0$ for all $x$).}
Finally, ${\Ima(A)}$ denotes the range of a matrix $A \in \mathbb{R}^{m \times n}$.

\section{Problem Formulation} \label{Problem_Formulation}

\subsection{Network model}\label{Network_Model}

 We describe the power network by a connected graph $(N,E)$, where $N = \{1,2,\dots,|N|\}$ is the set of buses and $E \subseteq N \times N$ the set of transmission lines connecting the buses.
We consider two types of buses in the network, buses with inertia and buses without inertia,
assuming non-trivial generation dynamics only in the first, since generators have inertia.
 We let $G = \{1,2,\dots,|G|\}$ and $L= \{|G|+1,\dots,|N|\}$ be the sets of buses with and without inertia respectively such that $|G| + |L| = |N|$.
Moreover, the term $(i,j)$ denotes the link connecting buses $i$ and $j$. The graph $(N,E)$ is assumed to be directed with an arbitrary direction, so that if $(i,j) \in E$ then $(j,i) \notin E$.
Additionally,  the sets of buses that precede and succeed bus~$j$ are denoted by $i:i\rightarrow j$ and $k:j\rightarrow k$ respectively.
 It should be noted that the form of the dynamics in~\eqref{sys1}--\eqref{dynsys} below is not affected by changes in graph ordering, and our results are independent of the choice of direction.  We make the following assumptions for the network: \newline
1) Bus voltage magnitudes are $|V_j| = 1$ p.u. for all $j \in N$. \newline
2) Lines $(i,j) \in E$ are lossless and {have susceptance with} 
{magnitude} $B_{ij} = B_{ji} > 0$. \newline
3) Reactive power flows do not affect bus voltage phase angles and frequencies. \newline

Following the above, we use the swing equation to describe the rate of change of frequency at generation buses. Moreover, power must be conserved at each of the load buses.
 This motivates the following system dynamics (e.g.~\cite{Bergen_Vittal}),
\begin{subequations} \label{sys1}
\begin{equation}
\dot{\eta}_{ij} = \omega_i - \omega_j, \; (i,j) \in E, \label{sys1a}
\end{equation}
\begin{equation}
M_j \dot{\omega}_j = - p_j ^L + p_j^M - \Lambda_j\omega_j - \sum_{k:j\rightarrow k} p_{jk} + \sum_{i:i\rightarrow j} p_{ij}, \; j\in G, \label{sys1b}
\end{equation}
\begin{equation}
 0 = - p_j ^L -  \Lambda_j\omega_j - \sum_{k:j\rightarrow k} p_{jk} + \sum_{i:i\rightarrow j} p_{ij}, \; j\in L, \label{sys1c}
\end{equation}
\begin{equation}
 p_{ij}=B_{ij} \sin\eta _{ij}, \; (i,j) \in E. \label{sys1d}
\end{equation}
\end{subequations}
In system~\eqref{sys1}, the time-dependent variables $\omega_j$ and $p^M_j$ represent respectively the deviation of the frequency at bus $j$ from its nominal value, namely 50Hz (or 60Hz), and
the mechanical power injection to the generation bus $j$.
 The positive constants $\Lambda_j$ and $M_j$ represent the frequency damping coefficient and generator inertia at any bus $j$ and generation bus $j$ respectively.
 The time-dependent variables $\eta_{ij}$ and $p_{ij}$ represent, respectively, the power angle difference\footnote{The  phase differences between buses $i$ and $j$ satisfy $\eta_{ij} = \theta_i - \theta_j$. The angles themselves must also satisfy $\dot{\theta}_j = \omega_j, j \in N$. This equation is omitted in \eqref{sys1} since the power transfers are functions of the phase differences only.}
 and the power transferred from bus $i$ to bus $j$.
 {Finally, $p_j^L$ denotes the uncontrollable demand at bus $j$.
 Below, we consider a wide class of generation dynamics and study {the} {stability properties of the equilibria of the system.}}

\begin{remark}
The analysis presented in this paper can be trivially extended to incorporate controllable demand. However, we focus on generation since the use of high order schemes is more relevant in this case and also for brevity in presentation.
\end{remark}

\subsection{Generation Dynamics}

To investigate control policies for a broad class of dynamics, we consider generation dynamics of the form
\begin{equation} \label{dynsys}
\begin{aligned}
&\dot{x}^M_j = A_jx^M_j + B_ju_j,\\
&p^M_j = C_jx^M_j + D_ju_j,
\end{aligned}, j \in G
\end{equation}
with input $u_j(t) \in \mathbb{R}$, state $x^M_j(t) \in \mathbb{R}^{n_j}$, $n_j \geq
 0, j \in G$,
output $p^M_j(t) \in \mathbb{R}$ and corresponding matrices $A_j \in \mathbb{R}^{n_j \times n_j}, B_j \in \mathbb{R}^{n_j}_j, C_j \in \mathbb{R}^{1 \times n_j}$ and $D_j \in \mathbb{R}$.
We assume in~\eqref{dynsys} that $A_j$ is Hurwitz which implies
that given any constant input {$u_j(t) = \bar{u}_j$}, there exists an
 asymptotically stable equilibrium point $\bar{x}^M_j \in \mathbb{R}^{n_j}$, such that $A_j\bar{x}^M_j + B_j\bar{u}_j = 0$.
 Correspondingly, there exists a constant
$K_j \in \mathbb{R}$, satisfying $K_j = -C_j A_j^{-1}B_j + D_j$, such that for any constant input $\bar{u}_j$ and corresponding state $\bar{x}_j$, the output $\bar{p}^M_j$ is given by
\begin{equation}\label{static_map_gen}
\bar{p}^M_j = C_j \bar{x}^M_j + D_j \bar{u}_j = K_j \bar{u}_j.
\end{equation}

Note that linear  {systems} are widely used in the literature to model generation {dynamics} (see e.g.  \cite[Section 11.1]{Bergen_Vittal}, \cite[Section 11.1.7]{machowski2011power}). Such models are particularly relevant when small disturbances are considered.
Furthermore, the Assumption that $A$ is Hurwitz, {i.e.,
generation dynamics are open-loop stable, 
is} in line with {practical implementations}.

The aim of this paper is to provide design conditions for the dynamics and control inputs of \eqref{dynsys} that ensure that secondary frequency control objectives are satisfied and stability and optimality are guaranteed.

\subsection{Optimal Generation Regulation}\label{sec: opt_problem}

We aim to study how generation should be adjusted to match  the {uncontrollable} demand  with minimum  cost.
Below we introduce an optimization problem, which we call the optimal generation regulation problem (OGR), that can be used to achieve this objective.

A quadratic cost function is used to describe the  cost induced when the generation output at bus $j$ is $p_j^M$, motivated from the fact that quadratic functions provide a local approximation for any convex cost function. Note that quadratic cost functions are commonly used in the literature \cite{zhao2015distributed, lavaei2012zero}.
The considered problem is to obtain the vector $p^M$ that minimizes the total cost and simultaneously satisfies a power balance constrain. In particular, we study the following optimization problem

\begin{equation} \label{Problem_To_Min}
\begin{aligned}
&\hspace{-2em}\underline{\text{OGR:}} \\
&\min_{p^M} \sum\limits_{j\in G} \frac{1}{2} q_j (p_{j}^M)^2,  \hspace{-1.5em}\\
&\text{subject to } \sum\limits_{j\in  G} p_j^M = \sum\limits_{j\in  N}  p_j^L,
\end{aligned}
\end{equation}
{where $q_j > 0, j \in N$ are the cost coefficients related with the generation cost at bus $j$.
{The equality constraint in~\eqref{Problem_To_Min} requires  the total generation, excluding some frequency dependent terms, to match the uncontrollable demand,
{which suffices to ensure that the frequency takes its nominal value at steady state.
The latter follows from summing \eqref{sys1b}--\eqref{sys1c} at equilibrium, which when combined with the equality constraint in~\eqref{Problem_To_Min} results to $\sum_{j \in N} \Lambda_j \omega_j^* = 0$, which allows to deduce that $\omega^*_j = 0, j \in N$ from \eqref{sys1a} at equilibrium and $\Lambda_j > 0, j \in N$.}
Note that more general quadratic cost functions can be considered, following similar approaches as in relevant literature, e.g. \cite{trip2017distributed}.
However, we consider the cost functions in \eqref{Problem_To_Min} for simplicity.}

\begin{remark}\label{rem_constraints}
The OGR problem could be extended by considering additional operational constraints, such as line capacities and fixed area power transfers.
{We chose not to consider this problem to keep the presentation compact and the focus of the paper on the problem presented in Section \ref{sec: Problem Statement} below. The investigation of such extensions is left as future work.} 
\end{remark}

\subsection{Problem Statement}\label{sec: Problem Statement}

In this section we present the problem we aim to solve.

\emph{Problem:}
Design the input $u_j, j \in G$ in \eqref{dynsys} such that it:
\begin{enumerate}[(i)]
\item Relies on local information.
\item Allows for distributed stability guarantees for \eqref{sys1}
when high order generation dynamics \eqref{dynsys}  are considered.
\item Ensures that the frequency takes its nominal value at steady state.
\item Ensures that system equilibria solve the OGR problem \eqref{Problem_To_Min}.
\item Relies on {easily obtainable} measurement requirements.
\item Is independent of (connected) network topology.
\end{enumerate}

{The first two
conditions {lead to a}
plug and play
operation of    high order generation dynamics, making the resulting scheme scalable
 and applicable in practise.}
 {In particular, the first condition requires that the control scheme uses information that can be locally measured at each bus or that can be transmitted from {neighbouring} buses.}
The third requirement is the main objective of secondary frequency control, i.e. to ensure that the frequency takes its  nominal value at equilibrium. Furthermore, condition (iv) requires optimality guarantees at steady state. The last two conditions  are related with the implementability of the designed controller, requiring it to {avoid hard to obtain measurement requirements, such as real time load measurements,} and be independent of  the topology of the network.

{As discussed in the introduction,} existing studies that attempted to resolve this problem {in a distributed fashion} had to either relax (ii) by assuming static generation dynamics (DAPI schemes), ignore (v) {imposing additional
measurement requirements} (Primal-Dual schemes), or relax (ii) by only considering first or second order generation dynamics and (v) by requiring measurements of generators internal states \cite{trip2017distributed}.
Our results are validated with simulations on more advanced dynamics than \eqref{sys1}--\eqref{dynsys}, as described in Section \ref{Simulation_NPCC}.

\section{Distributed averaging dynamic output controller}\label{sec:Power_Command_Dynamics}

In this section we present a novel scheme for {distributed optimal secondary frequency regulation,
which we will refer to as}
 distributed averaging dynamic output controller (DADOC). The proposed scheme offers advantages over existing distributed secondary frequency control schemes, {that will be} discussed in Section \ref{sec: Comparison}.

We consider a communication network described by a connected graph $(G, \tilde{E})$, where $\tilde{E}$ denotes the edges of the communication graph. The DADOC scheme proposed is given by equations \eqref{p_command}, \eqref{input_description} below. In particular
\begin{equation}\label{p_command}
\gamma_{j} \dot{p}^c_j = p^M_j - K_j u_j - k_{f,j} \omega_j + \hspace{-3mm}\sum_{i: (i,j) \in \tilde{E}}\hspace{-3mm} \alpha_{ij} (p^c_i - p^c_j),  j \in G,
\end{equation}
where $p^c_j$ is a power command signal, $\alpha_{ij} = \alpha_{ji}, \gamma_j$ and $k_{f,j}$ are positive constant gains of the controller and {$K_j$ follows from \eqref{static_map_gen}}.
{The generation input is described by}
\begin{equation}\label{input_description}
u_j = k_{c,j}p^c_j -k_{d,j}\omega_j, j \in G,
\end{equation}
where $k_{c,j}, k_{d,j}$ are positive {design} constants.

The control scheme in \eqref{p_command} contains three sets of terms.
The terms $k_{f,j} \omega_j $ and $\sum_{i: (i,j) \in \tilde{E}}\alpha_{ij} (p^c_i - p^c_j)$ respectively guarantee that the equilibria of the system satisfy the objective of secondary frequency control, i.e. that $\omega^* = \vect{0}_{|N|}$, and $p^{c,*}_i = p^{c,*}_j$ for all $i,j \in N$,  a  desired feature used to obtain an optimality interpretation of the steady state power allocation (see also Lemma \ref{eqbr_lemma} below).
The term $p^M_j - K_j u_j$  is zero at equilibrium, as follows from \eqref{dynsys}.
However, as shall be seen in the following sections, this term has a pivotal role in providing stability guarantees when high order generation dynamics of the form \eqref{dynsys} are considered.

We choose the generation input $u_j$ to be a weighted sum of frequency and power command, and let the weight coefficients  be design parameters.
Relative to the problem statement in Section \ref{sec: Problem Statement}, it follows that the control scheme in \eqref{p_command}--\eqref{input_description} satisfies (i) and (v) by having only local measurement requirements which are easy to obtain in practice (frequency and generation).
In the following section, we provide conditions for the choice of these design parameters such that
 convergence to the nominal frequency is achieved  {for arbitrary network topologies} {when high order generation dynamics are involved} while also taking optimality considerations into account, {hence satisfying the remaining objectives of the considered problem.}

\begin{remark}\label{rem_generation_subset}
An important aspect of DADOC schemes is that they do not need to be implemented on load buses. 
Furthermore, the employment of DADOC schemes on any non-empty subset of $G$ suffices to ensure that the equilibrium frequency is equal to the nominal (see also Lemma \ref{eqbr_lemma} {and its proof}).
\end{remark}

\subsection{Equilibrium analysis}\label{sec_eqlbr_analysis}

We now define an equilibrium
of system~\eqref{sys1}, \eqref{dynsys}, \eqref{p_command}, \eqref{input_description}.

\begin{definition} \label{eqbrdef}
The point $\beta^* = (\eta^*, \omega^*,x^{M,*}, p^{c,*})$ defines an equilibrium of the system~\eqref{sys1}, \eqref{dynsys}, \eqref{p_command}, \eqref{input_description} if all time derivatives of ~\eqref{sys1}, \eqref{dynsys}, \eqref{p_command}, \eqref{input_description} are equal to zero at this point.
\end{definition}

Throughout the paper, we assume the existence of some equilibrium to~\eqref{sys1}, \eqref{dynsys}, \eqref{p_command}, \eqref{input_description}, denoted by $\beta^* = (\eta^*, \omega^*, x^{M,*}, p^{c,*})$, as defined in Definition~\ref{eqbrdef}.
Furthermore, we use $(p^*, p^{M,*}, u^*)$ to represent the equilibrium values of respective quantities
in~\eqref{sys1}, \eqref{dynsys}, \eqref{p_command}, \eqref{input_description}.
{Note that the existence of an equilibrium is associated with the presence of sinusoids in the power transfer term \eqref{sys1d}, suggesting that power transfers are bounded at each bus and hence arbitrary mismatches between generation and demand cannot be tolerated.}
The study of the existence of equilibria is beyond the scope of this paper and the interested reader is referred to {e.g.} \cite{bolognani2015existence}, \cite{dorfler2013synchronization}.

Lemma \ref{eqbr_lemma} below, proven in the appendix, concerning the system~\eqref{sys1},~\eqref{dynsys}, \eqref{p_command},~\eqref{input_description} suggests that {the} frequency attains its nominal value at equilibrium. Furthermore, it shows  that power command variables synchronize at steady state, a property that, as discussed below, allows an optimality interpretation of the resulting equilibria.

 \begin{lemma}\label{eqbr_lemma}
Any equilibrium point $\beta^*$ given by Definition \ref{eqbrdef} satisfies $\omega^{*} = \vect{0}_{|N|}$ and $p^{c,*} \in {\Ima}(\vect{1}_{|G|})$.
 \end{lemma}

Furthermore, the following condition is imposed on the equilibrium values of power angle differences. This assumption is common in the power networks literature and can be considered as a security constraint.

\begin{assumption} \label{assum_angles}
$| \eta^*_{ij} | < \tfrac{\pi}{2}$ for all $(i,j) \in E$.
\end{assumption}

\subsection{Optimality analysis} \label{sec:optim}

Within the paper, we aim to provide conditions on generation dynamics, described by \eqref{dynsys}, such that convergence to an optimal point of \eqref{Problem_To_Min} is guaranteed.
Proposition \ref{optthm} below provides conditions on how the controller gains in \eqref{p_command}--\eqref{input_description} should be selected such that the equilibria of the system are solutions\footnote{Note that an equilibrium point is a solution to the OGR problem when at that point {the value of $p^M$ is optimal for \eqref{Problem_To_Min}.}} to the considered optimization problem \eqref{Problem_To_Min}.
 We will then show in the following section, how convergence to optimality can be achieved.

\begin{proposition} \label{optthm}
Consider equilibria of \eqref{sys1}, \eqref{dynsys}, \eqref{p_command}, \eqref{input_description},  characterized by Lemma \ref{eqbr_lemma}. Then, if the control {gains $k_{c,j}$} in \eqref{input_description} are chosen such that
\begin{equation}\label{opt_cond}
k_{c,j} = \frac{1}{q_j K_{j}}, j \in G,
\end{equation}
holds, then the equilibrium values $p^{M,*}$ are optimal for the OGR problem \eqref{Problem_To_Min}.
\end{proposition}

\begin{remark}
{Proposition} \ref{optthm} provides conditions on the choice of the design variable $k_{c,j}$ such that the equilibria of system \eqref{sys1}, \eqref{dynsys}, \eqref{p_command}, \eqref{input_description}, characterized by Lemma \ref{eqbr_lemma}, are solutions to the OGR problem \eqref{Problem_To_Min}.
It should be noted that design variables $k_{f,j}$ and $k_{d,j}$ do not appear in the optimization problem since they correspond to gains in frequency deviation that becomes zero at steady state.
However, their choice is important on the stability properties of the system, as described in the subsequent section.
\end{remark}

\section{Stability analysis}\label{sec:Convergence}

This section contains our main convergence results.
In particular, we provide conditions on how the gains in \eqref{p_command}--\eqref{input_description} should be selected and show that when those are satisfied, then convergence is guaranteed.

\subsection{Controller design conditions}

In this section, we impose a condition involving design constants $k_{c,j}, k_{d,j}$ and $k_{f,j}$, which is used in the convergence theorem presented in Section \ref{sec:Main_result} below. We then explain how this condition can be numerically tested in a computationally efficient way. The considered condition is presented below.

\begin{design condition}\label{assum_gains}
For each generation bus $j$, with dynamics described by \eqref{sys1}, \eqref{dynsys}, \eqref{p_command}, \eqref{input_description}, the controller parameters $k_{c,j}, k_{d,j}$ and $k_{f,j}$ are such {that
\begin{equation}\label{strictness_condition}
\begin{bmatrix}
-K_jk_{c,j} + D_j k_{c,j} & r^T_j \\
r_j & \hat{Q}_j \\
\end{bmatrix}  \leq 0,
\end{equation}
where $r_j = \begin{bmatrix}
\frac{k_{c,j}B_j^TP_j + C_j}{2}  & \frac{k_{f,j}-k_{d,j}K_j + D_j k_{d,j} - D_jk_{c,j}}{2}
\end{bmatrix}^T$, and}
\begin{equation*}\label{primary_LMI}
\hat{Q}_j = \begin{bmatrix}
\frac{P_jA_j + A_j^TP_j}{2} & \frac{k_{d,j}P_jB_j - C_j^T}{2} \\
\frac{k_{d,j}B^T_jP_j - C_j}{2}  & -\hat{\Lambda}_j - D_jk_{d,j}
\end{bmatrix}
\end{equation*}
  holds for some $P_j = P_j^T > 0$ and some $\hat{\Lambda}_j < \Lambda_j$.
\end{design condition}

Design condition \ref{assum_gains} is the main stability condition imposed on this paper, and is feasible for a broad class of linear systems, as discussed in Section \ref{sec:Discussion}.

{A necessary condition for \eqref{strictness_condition} to hold is $\hat{Q}_j \leq 0$. The latter {is} a sufficient stability condition for primary frequency control when linear generation dynamics are considered (see \cite[Sec. III-C]{kasis2018novel}).}
 Therefore, part of the stability conditions imposed for secondary frequency control can be seen as conditions for stability in primary frequency regulation. Furthermore, the condition requires knowledge of only a lower bound of the local damping coefficient $\Lambda_j$, allowing the analysis of systems {where the frequency damping is unknown but a lower bound is known}.

\begin{remark}
The inequality condition \eqref{strictness_condition} is an LMI with respect to the matrix $P_j$ and design parameter  $k_{f,j}$ and can hence be  verified in a computationally efficient way.
Note that the flexibility in choosing $P_j$ and $k_{f,j}$  in Design condition \ref{assum_gains} can be exploited to form various design optimization problems. One such problem would be to obtain the minimum frequency damping $\Lambda_j$ such that \eqref{strictness_condition} is satisfied, when particular generation dynamics are considered. Therefore, Design condition \ref{assum_gains} can be useful in system design.
\end{remark}

\begin{remark}\label{rem_passivity}
Design condition \ref{assum_gains}, in conjunction with the DADOC scheme dynamics,
is associated with the network independent Lyapunov function, which is used in the proof of  Theorem \ref{combined_thm} to obtain decentralized stability guarantees. 
The choice of the Lyapunov function is inspired from passivity, a notion that is widely used to provide scalable stability conditions in networks.
In particular, note that when the communication term $\sum_{i: (i,j) \in \tilde{E}}$ $\alpha_{ij} (p^c_i - p^c_j)$ is omitted in \eqref{p_command} {and Design condition \ref{assum_gains} holds}, then the system with input $(-\omega_j)$ and output $(p^M_j -\Lambda_j \omega_j)$ is an input strictly passive system. 
\end{remark}

\subsection{Main result}\label{sec:Main_result}

We now state our main result, demonstrating local convergence to an optimal point of \eqref{Problem_To_Min} where the frequency attains its nominal value.

\begin{theorem}\label{combined_thm}
{Consider an equilibrium of~\eqref{sys1}, \eqref{dynsys}, \eqref{p_command}, \eqref{input_description} such that Assumption~\ref{assum_angles} holds and let Design condition \ref{assum_gains} and \eqref{opt_cond} be satisfied.
Then, there exists an open neighborhood of initial conditions about this equilibrium such that the solutions of~\eqref{sys1}, \eqref{dynsys}, \eqref{p_command}, \eqref{input_description}  asymptotically converge to
{a global minimum of}
the OGR problem \eqref{Problem_To_Min} with $\omega^* = \vect{0}_{|N|}$.}
\end{theorem}

Theorem \ref{combined_thm} demonstrates local convergence to an optimal solution of the OGR problem \eqref{Problem_To_Min} that, following Lemma \ref{eqbr_lemma}, satisfies $\omega^* = \vect{0}_{|N|}$, {hence satisfying all the objectives of the problem statement in Section \ref{sec: Problem Statement}.}
The main conditions for stability are Assumption~\ref{assum_angles}, which is abundant in the power literature, and Design condition \ref{assum_gains}.
In the following section, we explain how Design condition \ref{assum_gains} applies to various generation schemes.

\begin{remark}\label{rem_local}
It should be noted that the local nature of the convergence result in Theorem \ref{combined_thm} follows from the presence of sinusoids in the power transfers \eqref{sys1d}.
In particular, if \eqref{sys1d} is linearized, then the set of equilibria of the resulting system would be globally attractive.
\end{remark}

\section{Discussion}\label{sec:Discussion}

In this section we discuss the applicability of our proposed scheme relative to existing schemes considered in the literature. Furthermore, we provide several examples of relevant generator models that fit within the considered analysis.

\subsection{Comparison with existing literature}\label{sec: Comparison}

As discussed in the introduction, the problem of addressing issues of  stability and optimality for secondary frequency control in a distributed way has been widely considered  in the literature in recent years. Most studies focused on two particular approaches to address this problem.
The first approach ensures that {the} frequency takes its nominal value at equilibrium by using integral action on {the} frequency.
This approach, resulting to a distributed averaging proportional integral (DAPI) controller  has been considered in many recent studies in the literature \cite{trip2016internal, dorfler, dorfler2014plug, simpson2013synchronization,  andreasson2014distributed}.
 DAPI schemes are simple to implement,  since further than exchanging a synchronizing variable, they only require knowledge of the local frequency which is easily obtainable.
  However, the stability and optimality results along this setting are limited to static generation models and quadratic cost functions.

{The second approach that has been used guarantees that {the} frequency takes its nominal value at steady state by ensuring that the total generation, excluding certain frequency damping terms, is equal to the total demand.
This approach leads to a power command scheme that follows from a primal/dual algorithm associated with some optimization problem \cite{mallada2017optimal}, \cite{kasis2017stability}, \cite{low2014distributed}, \cite{zhao2015distributed}, \cite{li2016connecting}.}
Primal-Dual schemes
allow for stability guarantees when high order and {nonlinear} generation dynamics are included
and for
 economic optimality to be deduced when general convex cost functions are considered.
 However, the implementation of such controllers requires real time measurements of demand, which can be difficult to obtain in many cases, and generation.
 {Attempts to alleviate the requirement for demand measurements introduced additional information flow requirements (power transfers, frequency derivatives) and {knowledge} of certain parameters (frequency damping, generation inertia).
 Such information requirements may introduce additional implementation costs and possibly raise security issues if the quality of information is compromised, e.g. due to  cyber-attacks.

DADOC schemes ensure that {the} frequency takes its nominal value at equilibrium by integrating the difference between the generation output and its static map for given power command input and zero frequency deviation.
Their structure enables them to share benefits of both mentioned control schemes.
In particular, DADOC schemes, {in contrast with DAPI schemes,
{allow for}  stability guarantees to be obtained when first or higher order generation dynamics are considered for a general network by tuning local control variables.}
Its only additional, but not restrictive, requirement compared to DAPI schemes is that of generation output measurements.
{Compared to DADOC schemes,} Primal-Dual schemes additionally allow to incorporate nonlinear generation dynamics and general convex cost functions.
Table \ref{table_comparison} summarizes the comparison between DADOC, DAPI and Primal-Dual schemes.

An additional feature of DADOC schemes is that the required equilibrium condition for frequency, i.e. that $\omega^* = \vect{0}_{|N|}$, is achieved without requiring {information from} all buses, which is important when a controller at a bus is withdrawn due to a failure {(see also Remark \ref{rem_generation_subset})}.

\begin{table}
\begin{center}
\resizebox{0.49\textwidth}{!}{
 \begin{tabular}
 {|c | c | c | c |}
 \hline
  & DAPI & DADOC & Primal-Dual \\ [0.5ex]
 \hline
 Allowable cost  & Quadratic & Quadratic & Convex \\
 function models & & & \\[0.5ex]
 \hline
 Allowable generation & Static & High Order & High Order \\
 dynamics & & & and {Nonlinear} \\[0.5ex]
 \hline
  {Information flow} & {Very low} & {Low}  & {High} \\
 {requirements} & &  &  \\[0.5ex]
 \hline
\end{tabular}
}
 \caption{Comparison between the two dominant schemes in {the} literature, DAPI and Primal-Dual, and DADOC schemes in terms of allowable cost function models and generation dynamics and measurement requirements.}
 \label{table_comparison}
 \end{center}
\end{table}

A further attempt to address distributed secondary frequency regulation issues has been made in \cite{trip2017distributed}, where the proposed controller  allows for stability guarantees to be obtained when first or second order generation dynamics are considered, using measurements of generators internal states.
DADOC schemes have  {milder} measurement requirements (generation output instead of internal states) for their implementation compared to the scheme proposed in  \cite{trip2017distributed} and provide design conditions for power networks with  {arbitrary (finite)} order generation dynamics.
Furthermore,  when second order schemes are considered,
the imposed conditions on \cite{trip2017distributed} can be conservative, being applicable only if the  time constants ratio is less than $\min(\Lambda_j,4)$.
 Our analysis {allows any time constant ratio}, as follows from Lemma \ref{suff_cond_lemma} below.

\subsection{Applications of main results}\label{sec:Applications}

To demonstrate the {relevance} of our analysis, we present examples of first, second and fifth order turbine governor dynamics and explain how our proposed conditions apply to them.

Consider the first order generation dynamics described by
\begin{equation}\label{first_order}
\tau_j \dot{p}^M_j = -p^M_j + K_j(k_{c,j}p^c_j - k_{d,j}\omega_j),
\end{equation}
for some constant $\tau_j > 0$, coupled with the controller \eqref{p_command} and some frequency damping $\Lambda_j$.
{For this system, the controller parameters can always be suitably selected such that   both Design condition \ref{assum_gains} and optimality condition \eqref{opt_cond} are satisfied, as follows from Lemma \ref{first_order_lemma} below, proven in the appendix.

\begin{lemma}\label{first_order_lemma}
Consider a bus $j$ with generation dynamics described by \eqref{first_order} coupled with the power command dynamics \eqref{p_command}. Then, for any positive values for $\tau_j,q_j, K_j$ and $\Lambda_j$
 there  exist positive constants $k_{c,j}, k_{d,j}$ and $k_{c,j}$ such that   both Design condition \ref{assum_gains} and optimality condition \eqref{opt_cond} are satisfied.
\end{lemma}
}

The proposed framework also applies to higher order generation dynamics,
such as the following second order model describing turbine governor dynamics (e.g. \cite{Bergen_Vittal}),
\begin{subequations}\label{second_order}
\begin{equation}
 \dot{\alpha}_j = -\frac{1}{\tau_{a,j}}(\alpha_j - K_j(k_{c,j}p^c_j - k_{d,j}\omega_j)),
\end{equation}
\begin{equation}
 \dot{p}^M_j = -\frac{1}{\tau_{p,j}}({p^M_j} - \alpha_j),
\end{equation}
\end{subequations}
where $\alpha_j$ and $\tau_{a,j}, \tau_{p,j} > 0$ are the internal state of the system and time constants associated with the generation dynamics respectively. We considered the case where \eqref{second_order} is coupled with the power command dynamics described by \eqref{p_command} and some frequency damping $\Lambda_j$.
{Lemma \ref{suff_cond_lemma} below, proven in the appendix, provides a sufficient condition for the  value of frequency damping $\Lambda_j$ such that Design condition \ref{assum_gains} holds for the considered system.

\begin{lemma}\label{suff_cond_lemma}
Consider a bus $j$ with generation dynamics described by \eqref{second_order} coupled with the power command dynamics \eqref{p_command} and some frequency damping $\Lambda_j$. Then, Design condition \ref{assum_gains} holds for any positive values of $\tau_{a,j}, \tau_{p,j}$ if  \[\Lambda_j > \frac{K_j}{3k_{c,j}}(k_{c,j}^2 - k_{c,j}k_{d,j} + k_{d,j}^2)\] holds. 
\end{lemma}

The above lemma provides a sufficient condition for the value of frequency damping $\Lambda_j$ such that Design condition \ref{assum_gains} holds at a bus with generation and power command dynamics described by \eqref{second_order} and \eqref{p_command} respectively.
Note that this condition can {be relaxed} when particular values for $\tau_{a,j}, \tau_{p,j}$ are considered {and that it does not impose any constraint on time constants $\tau_{a,j}, \tau_{p,j}$.}}

{
To further demonstrate the applicability of our proposed scheme, we consider the fifth order turbine governor dynamics provided by the Power System Toolbox \cite{cheung2009power}.
The dynamics relating turbine governor power output {${p}^M_j$}  with the negative frequency deviation {$-{\omega}_j$} are described by the following transfer function,
\[
G_j(s)=K_j\frac{1}{(1+sT_{s,j})}\frac{(1+sT_{3,j})}{(1+sT_{c,j})}\frac{(1+sT_{4,j})}{(1+sT_{5,j})},
\]
where $K_j$ and $T_{s,j}, T_{c,j}, T_{3,j}, T_{4,j}, T_{5,j}$ are the droop coefficient and time-constants respectively.
We have studied the implementation of our proposed scheme on such dynamics using
realistic values for these models provided by the toolbox for the NPCC network\footnote{The data were obtained from {the} Power System Toolbox~\cite{cheung2009power} data file datanp48.}, where turbine governor dynamics are implemented on 22 buses.
The corresponding buses also have appropriate frequency damping $\Lambda_j$.
We examined the effect of incorporating an input signal as follows from \eqref{input_description} in the above dynamics by considering\footnote{Note that $\hat{x}$ denotes the Laplace transform of {$x$}.} the turbine governor dynamics
\[
\hat{p}^M_j=  G_j \hat{u}_j, j \in N.
\]
 For appropriate choices of design constants $k_{c,j}, k_{d,j}$ and $k_{f,j}$ in \eqref{p_command}, we have numerically validated that Design condition \ref{assum_gains} was satisfied at all 22 buses with turbine governor dynamics,
 {hence demonstrating the applicability of our analysis to realistic high order turbine governor dynamics.}}

\section{Simulation on the NPCC 140-bus system} \label{Simulation_NPCC}

In this section we use the Power System Toolbox~\cite{cheung2009power} to perform numerical simulations on
the Northeast Power Coordinating Council (NPCC) 140-bus interconnection system, in order to numerically validate our analytic results.
The model used by the toolbox is more detailed than our analytic one, including  a DC12 exciter model, a transient reactance generator model, and {high} order turbine governor models\footnote{The simulation details can be found in the data file datanp48 and the Power System {Toolbox  manual \cite{cheung2009power}}.}.

The NPCC network  consists of 47 generation and 93 load buses and has a total real power of 28.55GW.
For our simulation, we considered a step increase in demand of magnitude {$4$} p.u. (base 100MVA) at each of the load buses {2, 9, 16, 17, 19 and 20}  at $t=1$ second.

To demonstrate the applicability of DADOC control schemes, the dynamics in \eqref{p_command}--\eqref{input_description} where implemented on 11 generators with third, fourth and fifth order turbine governor dynamics.
Furthermore, a quadratic cost function, penalising the deviation in generation output was considered at each contributing generation bus.
The cost coefficients where selected to be equal to $K^{-1}_j$, relating high cost coefficients with small droop gains, in consistency with relevant literature on optimal frequency regulation, e.g. \cite{kasis2016primary}.
Furthermore, it has been numerically verified that the  stability and optimality properties of the system, demonstrated below, are retained for a broad range of cost coefficient values.
Note that the choice of controller parameters was in agreement with Design condition  \ref{assum_gains} and optimality condition \eqref{opt_cond}.

\begin{figure}[t]
\centering
\includegraphics[trim = 0mm 0mm 10mm 0mm, scale = 0.6,clip=true]{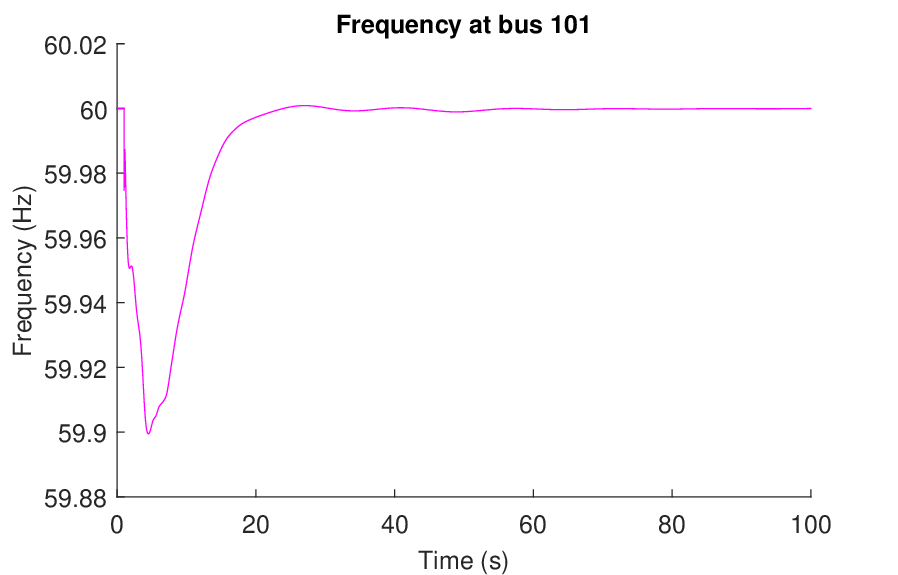}
\caption{Frequency at bus 101 of the NPCC network.}
\label{Frequency}
\end{figure}

\begin{figure}[t]
\centering
\includegraphics[trim = 0mm 0mm 10mm 0mm, scale = 0.6,clip=true]{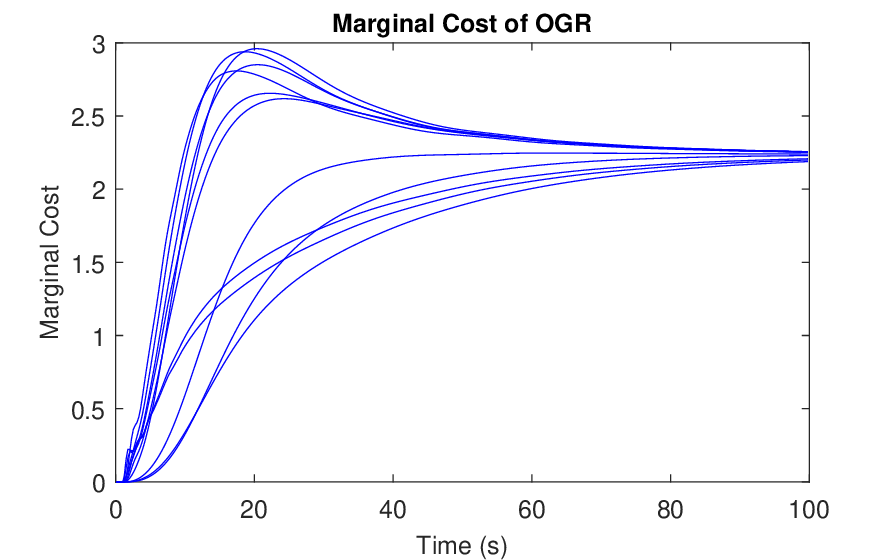}
\caption{Marginal cost of generation buses contributing to secondary frequency control.}
\label{Marginal}
\end{figure}

 Figure \ref{Frequency} depicts the frequency response on a randomly selected bus, where it is demonstrated that the frequency attains its nominal value at steady state.
  Hence, the simulation  numerically validates the analytic convergence results of Theorem \ref{combined_thm}.
Furthermore,  Figure \ref{Marginal} demonstrates
that the marginal costs
 of all 11 generators that contribute to secondary frequency control converge to the same value.
  This illustrates the optimality in the power allocation
among generators, since equality in the marginal cost is {sufficient} to solve \eqref{Problem_To_Min}.

{\subsection{Comparison of DADOC and DAPI schemes}\label{Simulations_comparison}

{DADOC schemes, as already discussed, provide distributed stability guarantees on power networks with high order generation dynamics. This is in contrast to DAPI schemes, the alternative scheme proposed in the literature with easily obtainable measurement requirements, that can only tolerate static generation.
 This subsection demonstrates the capability of DADOC schemes to offer improved performance compared to DAPI schemes.}

 {To compare the performance of DADOC and DAPI schemes, we repeated the simulation described above with $9$ generators performing secondary frequency control and introducing controllable loads described by \eqref{second_order} on 20 load buses. To excite the system, we considered a step disturbance of magnitude $2$ p.u. (base 100MVA) at $t=1$ second on buses 2, 9, 16, 17, 19 and 20.
 Furthermore, to make the simulations more realistic, we imposed  a communication delay of 170ms (see e.g. \cite{taylor2005wacs}) {on frequency measurements} and the transmission of power command signals  to generation units and controllers.
For a fair comparison, we aimed for the fastest achievable response from each scheme, {{(i.e, the} response with the minimum time for the} frequency to converge {to} within $0.01$Hz {of} its nominal value),  by accordingly adjusting the gains in the range where the system is well behaved.
The frequency response for both cases is depicted with the blue and green lines in Figure \ref{DADOC_DAPI_loads}, where
it is evident that DADOC schemes allow for a significantly faster response, since the frequency converges to within $0.01$Hz from its nominal value after $15.4$s compared to approximately $40$s when the DAPI scheme is implemented.
 {Furthermore, selecting the gains of DAPI scheme to provide an initial response of similar speed {(i.e.,} of similar time to reach the nominal frequency for the first time)  to DADOC} 
  resulted to oscillations, as also demonstrated in Figure~\ref{DADOC_DAPI_loads}.
}

\begin{figure}[t]
\centering
\includegraphics[trim = 0mm 0mm 9mm 0mm, scale = 0.7,clip=true]{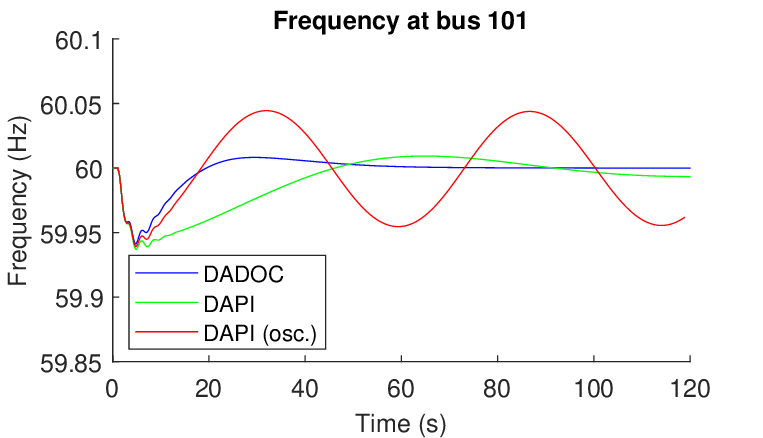}
\caption{{Frequency response at bus 101  when DADOC and DAPI schemes are implemented with fastest stable gains and time delays of 170ms (implementation with loads). The red line depicts the oscillating response following from applying similar gains to DAPI schemes as the ones tolerated by DADOC schemes.}}
\label{DADOC_DAPI_loads}
\end{figure}
}

\section{Conclusion}\label{Conclusions}

We have considered the problem of designing distributed control schemes for secondary frequency regulation such that stability is guaranteed and an optimal power allocation is attained.
We have presented a novel distributed averaging dynamic output control (DADOC) scheme that achieves the above objectives by guaranteeing stability and optimality when a wide class of generation dynamics and quadratic cost functions are considered, and that {the} frequency attains its nominal value at steady state.
DADOC controllers offer advantages compared to the existing schemes presented in the literature, allowing for stability guarantees when high order generation dynamics are incorporated and requiring easily obtainable measurements (generation output and frequency).
We demonstrate the relevance of our analysis  with several examples of turbine governor dynamics and  explain how design parameters can be selected in a computationally efficient way  by suitable LMI conditions.
Our analytic results are validated with realistic simulations on the NPCC 140-bus system, {where DADOC schemes demonstrate improved performance compared to DAPI schemes}.

  \renewcommand{\theequation}{A.\arabic{equation}}
  \setcounter{equation}{0}

\section*{Appendix}

This appendix contains the proofs of the results presented in this paper.

Throughout the proofs we will make use of the following equilibrium equations for the dynamics in \eqref{sys1}, \eqref{dynsys}, \eqref{p_command}, \eqref{input_description},
\begin{subequations}\label{eqbr}
\begin{gather}
0 = \omega^*_i - \omega^*_j, \; (i,j) \in E, \label{eqbr1} \\
0 = - p_j^L + p_j^{M,*} - \Lambda_j \omega^*_j - \sum_{k:j\rightarrow k} p^*_{jk} + \sum_{i:i\rightarrow j} p^*_{ij}, \; j\in G, \label{eqbr2} \\
0 = - p_j^L - \Lambda_j \omega^*_j - \sum_{k:j\rightarrow k} p^*_{jk} + \sum_{i:i\rightarrow j} p^*_{ij}, \; j\in L, \label{eqbr3} \\
p^{*}_{ij} = B_{ij} \sin \eta_{ij}^*, \; (i,j) \in E, \label{eqbr4} \\
p^{M,*}_j = K_{j} u^*_j, \; j \in G, \label{eqbr5} \\
0 = p^{M,*}_j - K_{j} u^*_j - k_{f,j} \omega^*_j  + \hspace{-3mm}\sum_{i: (i,j) \in \tilde{E}}\hspace{-3mm} \alpha_{ij} (p^{c,*}_i - p^{c,*}_j), j \in G, \label{eqbr6} \\
u^*_j =
k_{c,j} p^{c,*}_j -k_{d,j} \omega^*_j, j \in G. \label{eqbr7}
\end{gather}
\end{subequations}

\emph{Proof of Lemma~\ref{eqbr_lemma}:}
From \eqref{eqbr1} it follows that {$\omega^*_i = \omega^*_j$ for all $(i,j) \in E$, which results to $\omega^* \in \Ima(\vect{1}_{|N|})$}.
Then, summing \eqref{eqbr6} over all $j \in G$ results to
$
 \sum_{j \in G} (p^{M,*}_j -  K_ju^*_j - k_{f,j}\omega^*_j) = 0
$,
which by \eqref{eqbr5} and $k_{f,j} > 0$ implies that $\omega^* = \vect{0}_{|N|}$.
Since $\omega^*_j = 0, j \in N$, it follows by \eqref{eqbr5}, \eqref{eqbr6} and the fact that the communication graph is connected that $p^{c,*} \in {\Ima}(\vect{1}_{|G|})$.
\hfill $\blacksquare$

\emph{Proof of {Proposition}~\ref{optthm}:}
The OGR optimization problem~\eqref{Problem_To_Min} is convex and has a continuously differentiable cost function.
 Thus, a point $\bar{p}^M$ is a global minimum for~\eqref{Problem_To_Min} if and only if it satisfies the KKT conditions~\cite{boyd2004convex}
\begin{subequations} \label{kkt}
\begin{gather}
q_j \bar{p}^M_j = \nu, \; j \in G, \label{kkt1} \\
\sum\limits_{j\in  G} \bar{p}_j^M = \sum\limits_{j\in  N} p_j^L, \label{kkt2}
\end{gather}
\end{subequations}
for some constant $\nu \in \mathbb{R}$. It will be shown below that these conditions are satisfied by the equilibrium values $\bar{p}^M = p^{M,*}$ defined by equations~\eqref{eqbr5} and~\eqref{eqbr7} when \eqref{opt_cond} holds.

From Lemma \ref{eqbr_lemma}, it follows that $\omega^* = \vect{0}_{|N|}$ and $p^{c,*} \in {\Ima}(\vect{1}_{|G|})$.
Then, let $\nu = p^{c,*}_j$ and note that is common at every bus since power command variables synchronize at steady state.
 Therefore, it follows that $(q_j)^{-1} \nu  = (q_j)^{-1} p^{c,*}_j = (q_j K_j k_{c,j})^{-1}p^{M,*}_j = p^{M,*}_j$, by $\omega^* = \vect{0}_{|N|}$  and equations~\eqref{eqbr5}, \eqref{eqbr7} and~\eqref{opt_cond}. Thus, the optimality condition~\eqref{kkt1} holds.

Summing equations~\eqref{eqbr2} and~\eqref{eqbr3} over all $j \in G$ and $j \in L$ respectively and using that $\omega^* = \vect{0}_{|N|}$ shows that~\eqref{kkt2} also holds.
Hence, the values $\bar{p}^M = p^{M,*}$ satisfy the KKT conditions~\eqref{kkt}. Therefore, the equilibrium values $p^{M,*}$ define a global minimum for~\eqref{Problem_To_Min}.
\hfill $\blacksquare$

\emph{Proof of Theorem~\ref{combined_thm}:}
We will use the dynamics in \eqref{sys1}, \eqref{dynsys}, \eqref{p_command}, \eqref{input_description} and the matrices $P_j$ in Design condition \ref{assum_gains} to define a Lyapunov function for the system~\eqref{sys1}, \eqref{dynsys}, \eqref{p_command}, \eqref{input_description}.

Firstly, let $V_F (\omega^G) = \frac{1}{2}\sum_{j \in G} M_j (\omega_j - \omega^*_j)^2$. The time-derivative of $V_F$ along the trajectories of~\eqref{sys1}--\eqref{dynsys} is given by
\begin{align*}
\hspace{-0.5mm}\dot{V}_F &=\hspace{-1.5mm} \sum_{j \in N}\hspace{-0.5mm}(\omega_j - \omega^*_j) \Bigg(\hspace{-2.0mm}-p^L_j + p^M_j - \Lambda_j\omega_j -\hspace{-2.5mm} \sum_{k:j\rightarrow k} p_{jk} + \hspace{-2.0mm}\sum_{i:i\rightarrow j} p_{ij}\hspace{-0.5mm}\Bigg),
\end{align*}
by substituting~\eqref{sys1b} for $\dot{\omega}_j$ for $j \in G$ and adding extra terms for $j \in L$, which are equal to zero by~\eqref{sys1c}. Subtracting the product of $(\omega_j - \omega^*_j)$ with each term in~\eqref{eqbr2} and~\eqref{eqbr3}, this becomes
\begin{align}
\dot{V}_F =& \sum_{j \in G} (\omega_j - \omega^*_j) (p^M_j - p^{M,*}_j)
 - \sum_{j \in N}  \Lambda_j(\omega_j - \omega^*_j)^2 \nonumber \\
 +
&\sum_{(i,j) \in E} (p_{ij} - p^*_{ij}) (\omega_j - \omega_i), \label{VFdiff}
\end{align}
using the equilibrium condition~\eqref{eqbr1} for the final term.

Furthermore, let $V_C(p^c) = \frac{1}{2}\sum_{j \in N} \gamma_j(p^c_j - p^{c,*}_j)^2$. Using \eqref{p_command} the time derivative of $V_C$ can be written as
\begin{multline}
\dot{V}_C = \sum_{j \in N}  (p^c_j - p^{c,*}_j) \Big((p^M_j - p^{M,*}_j) -K_j(u_j - u^*_j) \\
 - k_{f,j}(\omega_j - \omega^*_j)
+\hspace{-3mm}\sum_{i: (i,j) \in \tilde{E}}\hspace{-3mm} \alpha_{ij} [ (p^c_i-p^{c,*}_i) - (p^c_j - p^{c,*}_j)] \Big).
\end{multline}

Additionally, define $V_P(\eta) = \sum_{(i,j) \in E} B_{ij} \int_{\eta^*_{ij}}^{\eta_{ij}} ( \sin \theta - {\color{black}\sin \eta^*_{ij}} ) \, d\theta$. Using~\eqref{sys1a} and~\eqref{sys1d}, the time-derivative is given~by
\begin{align}
\dot{V}_P &= \sum_{(i,j) \in E} B_{ij} (\sin \eta_{ij} - \sin \eta^*_{ij}) (\omega_i - \omega_j) \nonumber \\
&= \sum_{(i,j) \in E} (p_{ij} - p^*_{ij}) (\omega_i - \omega_j). \label{VPdiff}
\end{align}

Furthermore, from Design condition \ref{assum_gains},
it follows that there exist gains $k_{f,j},k_{c,j},k_{d,j}$ and a positive definite matrix $P_j = P^T_j$ such that \eqref{strictness_condition} holds. Then, {let $V^M_j(x^M_j) = \frac{1}{2}(x^M_j - x^{M,*}_j)^TP_j(x^M_j - x^{M,*}_j)$ and note that it is positive definite. Following \eqref{dynsys}, the time derivative of $V^M_j$ is given by
\begin{multline}\label{VM_derivative}
\dot{V}^M_j = \frac{1}{2}(x^M_j - x^{M,*}_j)^T(P_jA_j + A_j^TP_j) (x^M_j - x^{M,*}_j) \\
+ (x^M_j - x^{M,*}_j)^T B_j (u_j - u^*_j),
\end{multline}
where $u_j$ and $u^*_j$ are given by \eqref{input_description} and \eqref{eqbr7} respectively.}

Based on the above, we {consider the Lyapunov candidate
\begin{equation}\label{Lyapunov_function}
V(\eta, \omega^G, x^M, p^c) = V_F + V_P
 + \sum_{j \in G} V^M_j + V_C.
\end{equation}
Using \eqref{VFdiff}  - \eqref{VM_derivative}, the time derivative of V is given by
\begin{align}\label{dot_V}
\dot{V} =& \sum_{j \in G} \bigl[(\omega_j - \omega^*_j) (p^M_j - p^{M,*}_j)
 + \dot{V}^M_j  \nonumber\\
&\hspace{-7.5mm} + (p^c_j - p^{c,*}_j) \Big((p^M_j - p^{M,*}_j) -K_j(u_j - u^*_j)  - k_{f,j}(\omega_j - \omega^*_j)\Big)]
\nonumber \\
&\hspace{-7.5mm}
 - \sum_{j \in N} \Lambda_j(\omega_j - \omega^*_j)^2  - \sum_{(i,j) \in \tilde{E}} \alpha_{ij}((p^c_i - p^{c,*}_i) - (p^c_j - p^{c,*}_j))^2.
\end{align}
Using~\eqref{VM_derivative}, it therefore holds that
\begin{equation}\label{dot_V_2}
\dot{V} =  - \hspace{-2.5mm} \sum_{(i,j)\in \tilde{E}} \hspace{-2.0mm}  \alpha_{ij}((p^c_i - p^{c,*}_i) - (p^c_j - p^{c,*}_j))^2  + \sum_{j\in G}z_j^T Q_j z_j,
\end{equation}
where {$z_j = \begin{bmatrix}
(p^c_j - p^{c,*}_j) & (x^M_j - x^{M,*}_j)^T & -(\omega_j -\omega^*_j)
\end{bmatrix}^T$} and $Q_j$ is given by
\begin{equation}\label{Suff_cond_proof}
Q_j = \begin{bmatrix}
-K_jk_{c,j} + D_j k_{c,j} & r^T_j \\
r_j & \begin{smallmatrix}
\frac{P_jA_j + A_j^TP_j}{2} & \frac{k_{d,j}P_jB_j - C_j^T}{2} \\
\frac{k_{d,j}B^T_jP_j - C_j}{2}  & -\Lambda_j - D_jk_{d,j}
\end{smallmatrix}
\end{bmatrix}.
\end{equation}
Note that $Q_j$ {above} is identical to the matrix in \eqref{strictness_condition} when $\Lambda_j$ is replaced by some $\hat{\Lambda}_j < \Lambda_j$ and using Design condition \ref{assum_gains} on}
 \eqref{dot_V_2}, it therefore holds that
\begin{align}\label{vdotineq}
& \dot{V} \le - \hspace{-1.75mm} \sum_{j \in N} \bar{\Lambda}_j(\omega_j - \omega^*_j)^2  - \hspace{-2.75mm} \sum_{(i,j)\in \tilde{E}} \hspace{-2.25mm}  \alpha_{ij}((p^c_i - p^{c,*}_i) - (p^c_j - p^{c,*}_j))^2 \nonumber \\
 &\le 0
\end{align}
for $\bar{\Lambda}_j = \Lambda_j - \hat{\Lambda}_j > 0, j \in N$.
Furthermore, $\omega^{G,*}$ and  $x^{M,*}_j, j \in G$ are strict global minima for $V_F$ and  $V^M_j, j \in G$ respectively.
 Moreover,  $p^{c,*}$ is a strict global minimum of $V_C$.
 Furthermore, from Assumption~\ref{assum_angles} it follows that there exists some neighborhood of $\eta^*$ in which $V_P$ is increasing.
  Since the integrand is zero at the lower limit of the integration, $\eta^*_{ij}$, this immediately implies that $V_P$ has a strict local minimum at $\eta^*$. Hence $V$ has a strict local minimum  point at ${\Gamma^*} := (\eta^*, \omega^{G,*}, x^{M,*}, p^{c,*})$.
Therefore there exists a {connected} set $\Xi := \{(\eta,  \omega^G, x^M, p^c) \colon V \le \epsilon\}$ containing $\Gamma^*$ where, for sufficiently small $\epsilon > 0$, $V$ is a nonincreasing function of all the system states, as follows from ~\eqref{vdotineq}, and has a strict local minimum at ${\Gamma^*}$.  Therefore, $\Xi$ contains $\Gamma^*$ and is compact and positively invariant for~\eqref{sys1}, \eqref{dynsys}, \eqref{p_command}, \eqref{input_description}.

Lasalle's Invariance Principle can now be applied with the continuously differentiable function $V$ on the compact positively invariant set~$\Xi$.
This guarantees that all solutions of~\eqref{sys1}, \eqref{dynsys}, \eqref{p_command}, \eqref{input_description} with initial conditions $(\eta(0),  \omega^G(0), x^M(0), p^c(0)) \in \Xi$ converge to the largest invariant set within $\Xi \, \cap \, \{(\eta,  \omega^G, x^M, p^c) \colon \dot{V} = 0\}$.
We now consider this invariant set. If $\dot{V} = 0$ holds within $\Xi$, then~\eqref{vdotineq} holds with equality, hence  we must have $\omega = \omega^*$ for all $j \in N$ and $(p^c - p^{c,*}) \in {\Ima}(\vect{1}_{|G|})$.

 Then, summing \eqref{eqbr2}--\eqref{eqbr3} over all $j \in G$ and $j \in L$ respectively, it follows that $\sum_{j \in G} p^{M}_j =  \sum_{j \in N} (p^L_j + \Lambda_j \omega^*_j) = c_1$, where $c_1$ is constant.
Furthermore, summing \eqref{p_command} over all $j \in G$, it follows that $\sum_{j \in G} \dot{p}^c_j = \sum_{j \in G}  (p^M_j - K_j u_j - k_{f,j} \omega^*_j) = \sum_{j \in G}  (p^M_j - K_j k_{c,j}p^c_j + (K_jk_{d,j} - k_{f,j}) \omega^*_j) = c_2 -\sum_{j \in G}(K_j k_{c,j}p^c_j)$, for some constant $c_2$.
Letting $\hat{p}^c = \sum_{j \in G} p^c_j$, from $(p^c - p^{c,*}) \in {\Ima}(\vect{1}_{|G|})$ it holds that $\dot{\hat{p}}^c = c_2 - c_3 \hat{p}^c$, where $c_3 = \sum_{j \in G}(K_j k_{c,j})$. {Therefore, it is trivial to show that $\hat{p}^c$ converges to some constant value, and since $(p^c - p^{c,*}) \in \Ima(\vect{1}_{|G|})$, it follows that also $p^c$ converges to some constant value $\bar{p}^c$.
The convergence of $(\omega, p^c)$ to $(\omega^*, \bar{p}^{c})$ results to $u_j = \bar{u}_j$ and by the conditions imposed on \eqref{dynsys} to $x^{M}_j = \bar{x}^{M}_j$ for all $j \in G$, where $(\bar{u}_j, \bar{x}^M_j)$ are constants.
Finally,  the {fact that $(\omega, p^M)$ are constant} guarantees from \eqref{eqbr2}--\eqref{eqbr4} that $\eta$ and $p$ are also constant.  Furthermore, by summing \eqref{eqbr2}--\eqref{eqbr3} and using \eqref{eqbr7} and the synchronization of $p^c_j$ variables, it follows that $\bar{p}^c$ is unique and therefore equal to $p^{c,*}$, which also implies that $(\bar{u}_j, \bar{x}^M_j)$ are equal to $(u^*_j, x^{M*}_j)$.}
{This allows to conclude the convergence of } all solutions of~\eqref{sys1}, \eqref{dynsys}, \eqref{p_command}, \eqref{input_description} with initial conditions $(\eta(0), \omega^G(0), x^M(0), p^c (0))$ $\in \Xi$  to the set of equilibrium points defined in Definition~\ref{eqbrdef} and characterized by Lemma \ref{eqbr_lemma}.
 Finally, to conclude the convergence proof, we  choose for $S$ any open neighborhood of ${\Gamma^*}$ within $\Xi$.
  Finally,  as follows from {Proposition} \ref{optthm}, when \eqref{opt_cond} holds, then the described equilibria are solutions to \eqref{Problem_To_Min}.
\hfill $\blacksquare$

Proof of Lemma \ref{first_order_lemma}:
First note that when the selection $k_{c,j} = \frac{1}{q_j K_j}$ is made then optimality condition \ref{opt_cond} is satisfied.
Then, consider the matrix in \eqref{strictness_condition} with the matrices $A_j, B_j, C_j$ and $D_j$ following from \eqref{first_order} and let $k_{d,j} = k_{c,j}$, $k_{f,j} = K_j k_{c,j}$, $P = \frac{\tau_j}{K_j k_{c,j}} > 0$ and $0 = \hat{\Lambda} < \Lambda_j$.
Then, two of the eigenvalues of the resulting matrix are zero and the third is given by $-\frac{((K_jk_{c,j})^2 + 1)}{K_jk_{c,j}}$ which is negative for all positive values of $K_j, k_{c,j}$.
\hfill $\blacksquare$

Proof of Lemma \ref{suff_cond_lemma}:
The proof follows by analytically evaluating the eigenvalues of the matrix in \eqref{strictness_condition} with the matrices $A_j, B_j, C_j$ and $D_j$ following from \eqref{second_order} and\footnote{Note that these choices for $P_j$ and $k_{f,j}$ where numerically seen to minimize the required amount of frequency damping $\Lambda_j$ such that Design condition \ref{assum_gains} was satisfied when generation dynamics described by \eqref{second_order} were considered. Also, note that Design condition \ref{assum_gains} is also feasible for other choices of $P_j$ and $k_{f,j}$.}
selecting $P_j = \frac{1}{K_j k_{c,j}} \begin{bmatrix}
\tau_{a,j} & 0 \\
0  & \tau_{p,j}
\end{bmatrix}$ and  $k_{f,j} = \frac{K_j}{2}(k_{c,j} + k_{d,j})$. In particular three of the eigenvalues of the resulting matrix are {non-positive} and the fourth is negative for all positive values of $\tau_{a,j}, \tau_{p,j}$ if $\Lambda_j > \frac{K_j}{3k_{c,j}}(k_{c,j}^2 - k_{c,j}k_{d,j} + k_{d,j}^2)$ holds.
\hfill $\blacksquare$

\bibliographystyle{IEEEtran}
\bibliography{andreas_bib}

\end{document}